\documentclass[12pt]{article}
\usepackage{amsfonts,amsmath,graphicx}
\usepackage{epstopdf}
\usepackage[final]{epsfig}
\usepackage{float}
\usepackage{graphicx}
\usepackage{amsfonts}
\usepackage{xcolor}

\input pictexwd.tex

\usepackage{soul}
\setstcolor{red}

\newtheorem{theorem}{Theorem}[section]
\newtheorem{corollary}[theorem]{Corollary}
\newtheorem{proposition}[theorem]{Proposition}

\newcommand{\proofend}{$\Box$\bigskip}

\title{Differential geometry of space curves: forgotten chapters}

\author{Dmitry Fuchs, Ivan Izmestiev, Matteo Raffaelli, \\
Gudrun Szewieczek, Serge Tabachnikov}

\date{}

\begin{document}

\maketitle

\section{Introduction} \label{sect:intro}

Evolutes, involutes, and osculating circles of curves belong to the main notions of planar differential geometry, going back to Christiaan Huygens in the 17th century. They have many interesting properties, including the surprising Tait--Kneser theorem:  the osculating circles of a curve with monotonic curvature are nested -- see, e.g., \cite{BJT,GTT} or \cite[Chapter~10]{FTb}. 

In this article we are concerned with evolutes and involutes of space curves. Although much of the material presented is not new and can be found in classic books such as \cite{B,E,St}, we believe that a modern and unified treatment, complemented with several novel observations and results, may be useful. Besides, the theory is illustrated with the help of computer graphics, a tool obviously not available to the classical geometers.

%\textcolor{red}{While being an old topic,} 
The geometry of space curves remains highly relevant in modern mathematics. For example, it is closely related to the theory of completely integrable systems: the filament (a.k.a.\ binormal, smoke ring, local induction) equation is a completely integrable evolution of space curves, equivalent to the nonlinear Schr\"odinger equation, see, e.g., \cite{La}. 
Another important application is the study of curved origami, which can be informally described as folding paper along curves (as opposed to straight lines); see, e.g., \cite{DO,FT1}.

In this article we consider three kinds of evolutes and involutes. They were all studied in the early days of differential geometry, however the terminology is not canonical. Following \cite{F} and \cite{FT}, we define the evolute as the locus of centers of osculating spheres or, equivalently, the curve whose osculating planes are the normal planes of the given curve (in \cite{B,E} the term {\it evolute} means something else, and Uribe-Vargas, in his detailed study of the evolute in our sense \cite{UV}, prefers the term {\it focal curve}). 

Evolutes and involutes of space curves possess many familiar properties of evolutes and involutes of plane curves, but also have some unexpected features. We survey these properties and surprises and discuss a natural modification of the construction in which the normal planes are replaced by the rectifying planes. For the resulting curve we use a somewhat awkward term {\sl pseudo-evolute}. Properties of pseudo-evolutes offer further surprises. 

In the plane, involutes are constructed by wrapping a non-stretchable string around a curve. The same construction in space provides yet another definition of involute, and one defines the evolute of a curve as the result of the converse operation. We call this third version {\it Monge evolutes} and {\it Monge involutes} -- they were studied by Gaspard Monge in the late 18th century. 

One of our goals is to describe the intricate interrelations between these three kinds of evolutes. 

In what follows we present some -- but not all -- calculations behind the geometrical statements; they involve only elementary calculus, but in some cases are cumbersome. Readers are encouraged to perform the missing calculations on their own.

%We restrict ourselves to the three-dimensional case, although there are immediate generalizations to higher dimensions.

\section{Textbook material} \label{textbook}

The following facts are undoubtedly known to the majority of the readers, but we prefer to provide a brief survey of them to establish the settings, terminology, and notations.

\paragraph{Frenet apparatus.}
By a {\it curve} we  mean an arclength-parametrized, generic (this notion will vary, but will always include the absence of inflexion points), and smooth (${\mathcal C}^\infty$, possibly with cusps) curve $\xi=\xi(t)$, or $x=x(t),y=y(t), z=z(t)$ in ${\mathbb R}^3$. If a curve has cusps, then the arcs between the cusps are unit speed parametrized.

The {\it Frenet frame} $({\bf t}={\bf t}(t),{\bf n}={\bf n}(t),{\bf b}={\bf b}(t))$ consists of the {\it tangent vector} ${\bf t}=\xi'$, the {\it (principal) normal}, which is the unit vector in the direction of $\xi''$, and the {\it binormal} $\bf b=t\times n$. The dependence of $\bf t,n,b$ on $t$ is described by the {\it Frenet formulas} 
$$
\begin{array} {rlrl} {\bf t}'=&&k{\bf n}\\ {\bf n}'=&-k{\bf t}&&+\tau{\bf b}\\ {\bf b}'=&&-\tau{\bf n}\end{array}
$$
where $k$ and $\tau$ are, respectively, the {\it curvature} and the {\it torsion} of the curve $\xi$. We denote the {\it radius of curvature} $1/k$ by $r$.

There are well-known formulas for the curvature and torsion:
\begin{equation*}
\label{eqn:CurvTorsUnit}
k=\|\xi''\|, \quad \tau=\frac{\det(\xi', \xi'', \xi''')}{\|\xi''\|^2}.
\end{equation*}
If the parametrization is not unit-speed, then the formulas become slightly more complicated:
 \begin{equation}
\label{eqn:CurvTors}
k = \frac{\|\xi' \times \xi''\|}{\|\xi'\|^3}, \quad
\tau = \frac{\det(\xi', \xi'', \xi''')}{\|\xi' \times \xi''\|^2}.
\end{equation}

The plane spanned by $\bf t$ and $\bf n$ is called the {\it osculating} plane of the curve $\xi$, the plane spanned by $\bf n$ and $\bf b$ the {\it normal} plane of $\xi$, and the plane spanned by $\bf t$ and $\bf b$ the {\it rectifying} plane of $\xi$.

\paragraph{Developable surfaces.}

A surface in space is called {\it developable} if it is locally isometric to a plane. An informal description of a developable surface is that it is a surface that can be made by bending, but not folding, a piece of paper. Certainly, we have in mind the ``ideal paper", which is ideally bendable (without any bounds on the curvature), non-compressible and non-stretchable, i.e., the length of any curve drawn on the paper remains unchanged in the process of bending.

The theory of developable surfaces was developed (pardon the unintended pun) in the late 1700's by Euler and Monge. Let us translate the main results into modern language. 

A generic (in particular, nowhere planar) developable surface is \textit{ruled}: every point belongs to a unique straight line, which is fully contained in the surface. Moreover, the tangent plane is constant along each line; this property distinguishes developable surfaces in the class of all ruled surfaces. 

Generically, these lines, called \textit{rulings}, are tangent to a certain curve, called the \textit{regression edge} of the surface. There the surface is not smooth: all sections by planes, transverse to the regression edge, have cusps.

The latter property (also a characteristic one) provides a universal method for constructing developable surfaces: take an arbitrary (generic) curve $\xi$ in space, possibly with cusps, and consider the union of all its tangent lines. This union is a developable surface with regression edge $\xi$, and all (generic) developable surfaces can be obtained in this way. This surface is called the {\it tangent developable} of $\xi$.

There are two types of non-generic developable surfaces: cylinders and cones over an arbitrary curve.
When the regression edge has a cusp, the developable surface becomes approximately conical; when the regression edge escapes to infinity, the developable becomes approximately cylindrical.

Another universal construction of developable surfaces uses an arbitrary (generic) one-parameter family of planes: for such a family there exists a unique developable surface, tangent to all  the planes, and called the {\it envelope} of the family. These planes are the osculating planes of the regression edge. For each plane, tangency occurs along a whole line, which is the tangent of the regression edge.

The last construction provides a convenient analytic description of everything mentioned above. Namely, if $F(x,y,z;t)=0$ is the equation of the family of planes (with parameter $t$), then we can form three systems of equations:
\begin{equation}
\label{eqn:Envelope}
\left\{\hskip-6pt\begin{array} {rl} F(x,y,z;t)&=0,\\ F'(x,y,z;t)&=0;\end{array}\right.\left\{\hskip-6pt\begin{array}{rl} F(x,y,z;t)&=0,\\ F'(x,y,z;t)&=0,\\ F''(x,y,z;t)&=0;\end{array}\right.\left\{\hskip-6pt\begin{array}{rl} F(x,y,z;t)&=0,\\ F'(x,y,z;t)&=0,\\  F''(x,y,z;t)&=0,\\ F'''(x,y,z;t)&=0.\end{array}\right.
\end{equation}
If we exclude $t$ from the first system, then we get the equation of the surface. If we solve the second system with respect to $x,y,z$, then we obtain a parametric equation of the regression edge. The solutions of the last system provide coordinates of the cuspidal points on the regression edge. 

Notice that all  this may be repeated for the case when $F(x,y,z;t)=0$ is a family of surfaces, not necessarily planes, but we will not need that level of generality. 

\section{Evolutes}\label{evolutes}

Let $\xi=\xi(t)$ be a (generic) curve in space. There arise three families of  planes: the osculating planes, the normal planes, and the rectifying planes. Each of them has an envelope (the regression edge of the developable surface tangent to the family). The first of these three cases is not interesting: the envelope is the curve $\xi$ itself.

The envelope of the family of normal planes is the {\it normal developable}, and its regression edge 
is what we call the {\it evolute} of the curve: this is  similar to the definition of the evolute of a planar curve as the envelope of the family of its normals.

\paragraph{Equation of the evolute.}

Let $\xi=\xi(t)=(x(t),y(t),z(t))$ be an arclength-parametrized curve, and $e=e(t)$ be its  evolute.

\begin{proposition}
\label{e=}
$e=\xi+r{\bf n}+\displaystyle\frac{r'}{\tau}{\bf b}.$
\end{proposition}

{\bf Proof.} The proof is based on  formulas (\ref{eqn:Envelope}). We denote the coordinates of points  in ${\mathbb R}^3$ as $P=(X,Y,Z)$, and we use the dot product to make the formulas more compact. 

The family of normal planes to $\xi$ is described by the formula
$$
F(P;t)=\xi'(t) \cdot (P-\xi(t))=0.
$$
The derivatives of $F$ are
$$
F'(P;t)=\xi''(t) \cdot (P-\xi(t))-1,\quad F''(P;t)=\xi'''(t) \cdot (P-\xi(t))
$$
(we used the equalities $\xi'(t)\cdot \xi'(t)=1$ and $\xi''(t)\cdot \xi'(t)=0$). 

The resulting (middle) system of equations (2) becomes
\begin{align*}
\begin{cases}
\xi'(t) \cdot (P - \xi(t)) &= 0,\\
\xi''(t) \cdot (P - \xi(t)) &= 1,\\
\xi'''(t) \cdot (P - \xi(t)) &= 0.
\end{cases}
\end{align*}
Solving this system for $P$ gives the evolute $e(t)$.

The Frenet formulas imply
$$
\xi' = {\bf t},\ 
\xi'' = k {\bf n},\
\xi''' = -k^2 {\bf t} + k' {\bf n} + k \tau {\bf b},
$$
hence 
\[
P - \xi = \frac{1}{k} {\bf n} - \frac{k'}{k^2} \frac{1}{\tau} {\bf b},
\]
and finally
\[
e = \xi + \frac{1}{k} {\bf n} - \frac{k'}{k^2} \frac{1}{\tau} {\bf b} = \xi + r {\bf n} + \frac{r'}{\tau} {\bf b},
\]
as required 
 \proofend

In this calculation, we assume that $\tau\ne 0$, otherwise the evolute escapes to infinity.

\paragraph{Osculating circles and osculating spheres.}

Let $\xi=\xi(t)$ be a curve, and let $\xi_0=\xi(t_0)$ be a point on this curve. For $t_1<t_2<t_3$ close to $t_0$ we denote by $C_{t_1,t_2,t_3}$ a circle passing through $\xi(t_1),\xi(t_2),\xi(t_3)$; for $t_1<t_2<t_3<t_4$ close to $t_0$ we denote by $S_{t_1,t_2,t_3,t_4}$ a sphere passing through $\xi(t_1),\xi(t_2),\xi(t_3),\xi(t_4)$. 

In the generic case, both $C_{t_1,t_2,t_3}$ and $S_{t_1,t_2,t_3,t_4}$ are well defined. Moreover, both have limits when all $t_i\to t_0$. This limit circle and limit sphere are called the {\it osculating circle} and the {\it osculating sphere} of the curve $\xi$ at the point $\xi_0$. 

It is well known (and obvious) that the osculating sphere can be described as the unique sphere which has a tangency of order $\ge3$ with $\xi$ at $\xi_0$. The osculating circle is the intersection circle of the osculating sphere with the osculating plane; it has tangency of order $\ge2$ with $\xi$ at the point $\xi_0$.

For our purposes, a more convenient description of the (center of the) osculating sphere is the following. We take three parameter values, $t_1<t_2<t_3$ close to $t_0$ and consider the three normal planes to the curve $\xi$ at points $\xi(t_1),\xi(t_2),\xi(t_3)$. These three planes have a common point, and  this point approaches the center of the osculating sphere at the point $\xi_0=\xi(t_0)$ when $t_1,t_2,$ and $t_3$ approach $t_0$.

\paragraph{Evolutes and osculating spheres.}

\begin{proposition}\label{locus} 
The  evolute of a curve is the locus of the centers of its osculating spheres.
\end{proposition}

{\bf Proof.} Let $F(X,Y,Z;t)=0$ be the equation of the normal plane to the curve $\xi$ at  point $\xi(t)$. Then the center of the osculating sphere at $\xi(t)$ is the limit for $\varepsilon\to0$ of the solution of the system
$$\left\{\begin{array} {ll} F(X,Y,Z;t-\varepsilon)&=0,\\ F(X,Y,Z;t)&=0,\\ F(X,Y,Z;t+\varepsilon)&=0,\end{array}\right.$$
which is equivalent to the system
$$\left\{\begin{array} {ll} F(X,Y,Z;t)&=0,\\ \displaystyle\frac{F(X,Y,Z;t+\varepsilon)-F(X,Y,Z;t-\varepsilon)}{2\varepsilon}&=0,\\ \displaystyle\frac{F(X,Y,Z;t+\varepsilon)-2F(X,Y,Z;t)+F(X,Y,Z;t-\varepsilon)}{2\varepsilon^2}&=0.\end{array}\right.$$ 
In particular, in the limit  $\varepsilon\to0$, this becomes
$$
\left\{\begin{array} {ll} F(X,Y,Z;t)&=0,\\ F'(X,Y,Z;t)&=0,\\ F''(X,Y,Z;t)&=0.\end{array}\right.
$$
The solution of the last system is the point $e(t)$ of the evolute. Thus $e(t)$ is the center of the osculating sphere of $\xi$ at $\xi(t)$. \proofend

Propositions \ref{e=} and \ref{locus} imply

\begin{corollary}\label{Radius}
Let $R$ be the radius of the osculating sphere of $\xi$ at $\xi(t)$. Then
$$
R^2=r^2+\left(\frac{r'}\tau\right)^2.
$$
\end{corollary}

\paragraph{Singularities of the  evolute.}

From Proposition \ref{e=} one has
$$
e'={\bf t}+ r'{\bf n}-kr{\bf t}+\tau r{\bf b}+\left(\frac{r'}{\tau}\right)'{\bf b}-\frac{r'}\tau\tau{\bf n}=\left(r\tau+\left(\frac{ r'}\tau\right)'\right){\bf b}.
$$
Thus, analogously to the two-dimensional case, the tangent to the evolute is always parallel to the binormal. It follows that  the cusps of the evolute occur when 
$$
\sigma:=r\tau+\displaystyle\left(\frac{r'}\tau\right)'=0.
$$
In particular, $\sigma \equiv 0$ is the condition for a curve to be spherical; indeed, the evolute of a spherical curve is a point, the center of the sphere.

On the other hand, from Corollary \ref{Radius},
$$
(R^2)'=2r'r+2\frac{ r'}\tau\left(\frac{ r'}\tau\right)',
$$
and so $(R^2)'=2\frac{r'}{\tau}\sigma$. This shows that at every cusp of the evolute, $R^2$ has zero derivative; thus, generically, $R$ achieves a maximum or a minimum.

But there is also the possibility that $r'=0$ but $\sigma\ne0$; so, although $R$ is maximal or minimal, the evolute has no cusp at this point. By Corollary \ref{Radius}, in this case $r=R$, that is, the osculating circle is a great circle  of the osculating sphere, and the center of the osculating sphere is contained in the osculating plane of the curve. This is one of the essential differences between the evolutes of planar and space curves.
\smallskip

 {\bf Remark.}
The quantity
\[
\frac{k^3\tau^2\sigma}{R^{5/2}}
\]
is called the {\it conformal torsion} and is one of the two invariants of space curves in conformal geometry (the other one being the conformal curvature); see, e.g., \cite{CSW}. In conformal geometry, spheres are ``flat", and the conformal torsion measures the deviation of the curve from its osculating sphere.

\paragraph{Examples.}

\begin{figure}[hbtp]
\centering
\includegraphics[width=.6\textwidth]{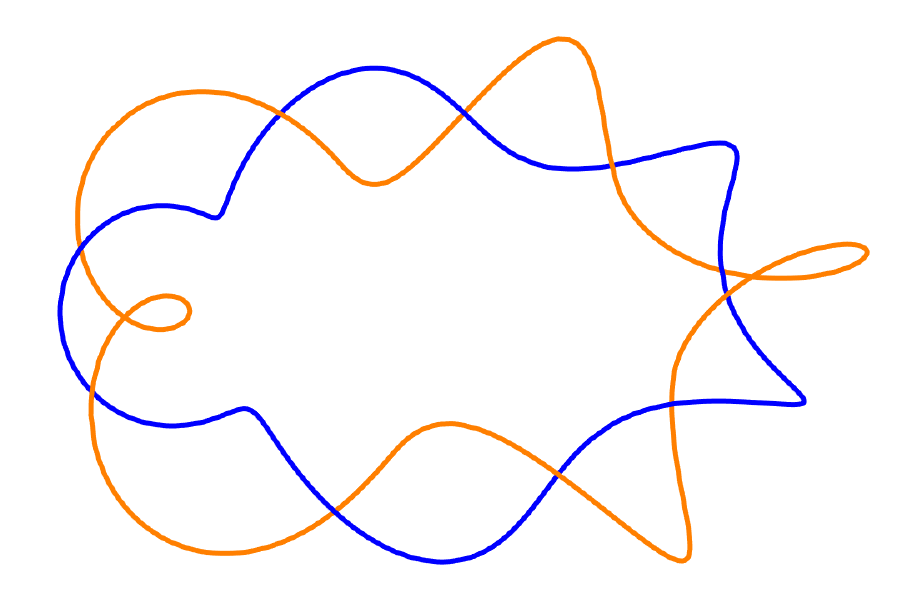}
\caption{A cusp-free evolute (orange) of a closed space curve (blue). It is a $(1,5)$-torus knot on a torus with minor radius $m= 0.15$ and major radius $l=1$. It is given by the parametrization $t \mapsto ((l+m \cos 5t)\cos t, (l+m \cos 5t)\sin t, -m \sin 5t )$.}
\label{torus}
\end{figure}

{\it A closed curve whose evolute has no cusps.} The evolute of a closed convex plane curve has cusps, in fact, at least four of them, according to the 4-vertex theorem (see, e.g., \cite{GTT} or \cite[chapter 10]{FTb}). But the evolute of a closed space curve with non-vanishing curvature and torsion may be free of cusps, see Figure \ref{torus}. 
\smallskip

\begin{figure}[hbtp]
\centering
\includegraphics[width=.5\textwidth]{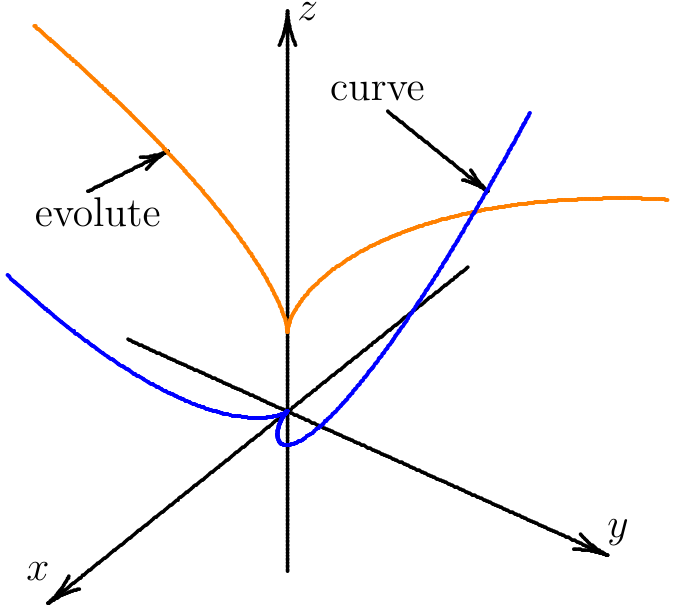}\qquad \quad
\includegraphics[width=.4\textwidth]{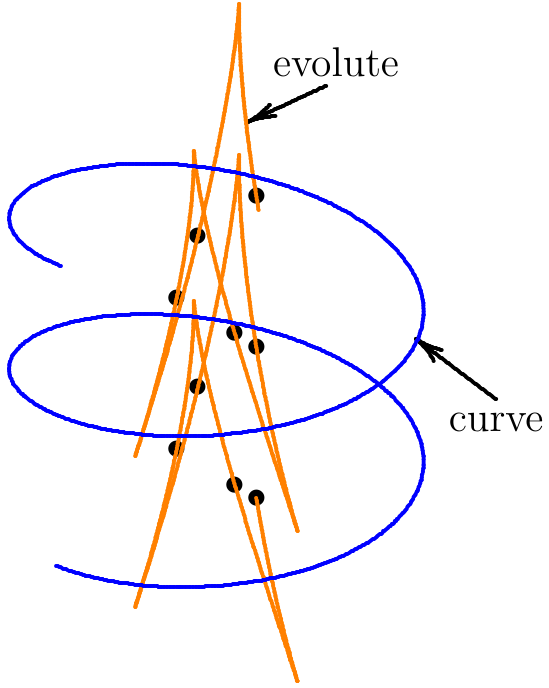}
\caption{Left: the evolute of a curve with a cusp. Right: the evolute of an elliptical helix.}
\label{ex23}
\end{figure}

{\it Evolute of a curve with a cusp.} It is well known that for a planar curve with a cusp (like a semi-cubic parabola), its evolute passes through the cusp and has no cusp at this point. The situation for a space curve is entirely different. Here is a parametric equation of a curve with a cusp and its evolute (see  Figure \ref{ex23} left):
$$
{\rm curve:}\quad\begin{array} {rl} x=&t^2,\\ y=&t^3, \\ z=&t^4;\end{array}\quad{\rm evolute:}\quad\begin{array} {rl} x=&\displaystyle\frac92t^4+20t^6,\\ y=&-8t^3-32t^5,\\ z=&\displaystyle\frac12+\displaystyle\frac92t^2+15t^4.\end{array}$$
We see that our curve has a cusp at the point $(0,0,0)$, and the evolute has a cusp at the point $(1/2,0,0)$.
\smallskip

{\it An elliptical helix.} The evolute of the standard (circular) helix is just another helix. For the {\it elliptical helix} $x=a\cos t,\, y=b\sin t, z=ct$, the evolute has 4 cusps for each turn (see Figure \ref{ex23} right). The heavy dots on the evolute mark its non-cuspidal points with $r'=0$ and $\sigma \ne 0$, which correspond to maxima or minima of $R$.

\paragraph{Interior and exterior points of a curve.}

Since at a generic point a curve has an odd degree tangency with the osculating sphere, a neighborhood of such point is contained either in the interior or the exterior of this sphere. It remains unclear how to visualize the difference between ``interior" and ``exterior" points. Just imagine that you have a rigid curve, say, a twisted bicycle spoke. Can you tell the interior points from the exterior ones? 

Once again let us use  system (2). Assuming that the three middle equations hold, the type of the point, that is, whether it is interior or exterior, is determined by the sign of $F'''(x,y,z;t)$. A calculation that we do not present yields the following result: a point is interior if and only if 
$$
\frac{\det({\bf t},{\bf t''},{\bf t}''')}{k^2\tau} > k^2.
$$

The points where a curve has a higher-order contact with its osculating sphere correspond to the singularities of the  evolute, that is, to the points where $\sigma = 0$.
This suggests that the above inequality can be rewritten in terms of $\sigma$.
Indeed, it turns out that $\det({\bf t},{\bf t''},{\bf t}''') = k^3 \tau^2 \sigma + k^4\tau$, so that the above inequality is equivalent to $\sigma \tau > 0.$

And here is a geometrical interpretation: a curve of positive torsion locally lies inside its osculating sphere if the center of the osculating sphere moves in the direction of the binormal to the curve, and the curve lies outside of its osculating sphere if the center of the sphere moves in the direction opposite to the binormal.

\paragraph{Curvature and torsion of the  evolute.} As we know, $e'= \sigma {\bf b}$.
Using the Frenet formulas, we calculate the next two derivatives of $e$:
$$
e''=\sigma'{\bf b}-\sigma\tau{\bf n},\quad e'''=\sigma k\tau{\bf t}-(2\sigma'\tau+\sigma\tau'){\bf n}+(\sigma''-\sigma\tau^2){\bf b}.
$$
The formulas for $e'$ and $e''$ show that the Frenet frame of the evolute has the form $({\bf t}_e, {\bf n}_e, {\bf b}_e) = (\pm {\bf b}, \pm {\bf n}, \pm {\bf t})$.

Formulas \eqref{eqn:CurvTors} imply

\begin{proposition}
\label{prp:CurvTorsEvol}
The curvature and the torsion of the  evolute are related to the curvature and the torsion of the initial curve by the formulas
\[
k_e = \frac{|\tau|}{|\sigma|}, \quad \tau_e = \frac{k}{\sigma}, \quad \text{where } \sigma=r\tau+\left(\frac{ r'}\tau\right)' \text{ and } r = \frac{1}{k}.
\]
\end{proposition}

Assume that the evolute has no cusps, that is, $\sigma \ne 0$. Without loss of generality, $\sigma > 0$. The magic cancellations when integrating the curvature or the torsion of the evolute
\begin{gather*}
\int_a^b k_e \|e'\|\, dt = \int_a^b \frac{|\tau|}{|\sigma|} |\sigma|\, dt = \int_a^b |\tau|\, dt,\\
\int_a^b \tau_e \|e'\|\, dt = \int_a^b \frac{k}{\sigma} |\sigma|\, dt = \int_a^b \mathrm{sgn}(\sigma) k\, dt
\end{gather*}
imply

\begin{corollary}
\label{cor:TotalCurvTors}
The total curvature of the evolute is equal to the total absolute torsion of the curve.
The total torsion of the evolute is equal to the total curvature of the curve taken with the sign of $\sigma$.
\end{corollary}

\paragraph{Curves congruent to their evolutes.} 
The planar curves congruent to their evolutes (with the congruence compatible with the parametrization) are logarithmic spirals. What about space curves?

\begin{proposition}
The only space curves congruent to their  evolutes (with the congruence compatible with the parametrization) are helices of equal curvature and torsion, and the curves whose curvature and torsion depend on an arclength parameter in the following way:
\[
k = \tau = c t^{-1/2}, \quad c > 0,
\]
as well as their mirror images.
\end{proposition}

\textbf{Proof.}
A necessary condition for a congruence compatible with the parametrization is $|\sigma| = 1$. 

Note that the reflection in a plane preserves $k$, but changes the sign of $\tau$, therefore it changes the sign of $\sigma$. Thus, without loss of generality, we may assume that $\sigma = 1$.

A space curve is determined up to congruence by its curvature and torsion as functions of an arclength parameter.
The formulas in Proposition \ref{prp:CurvTorsEvol} imply that $k = \tau$ and $\sigma = 1$ is a necessary and sufficient condition for a curve to be congruent to its evolute.

Substituting $k = \tau$ into the formula for $\sigma$, we obtain
\[
\sigma = 1 + (r'r)' = 1 + \left(\frac{r^2}{2}\right)''.
\]
Thus $\sigma = 1$ if and only if $r^2 = at + b$.
If $a = 0$, then $k = \tau = \mathrm{const}$, and one obtains helices of equal curvature and torsion.
If $a \ne 0$ then, after a time shift and, if necessary, time reversal, one has $r^2 = t/c^2$, that is, $k = \tau = c t^{-1/2}$.
\proofend

\begin{figure}[hbtp]
\centering
\includegraphics[width=.8\textwidth]{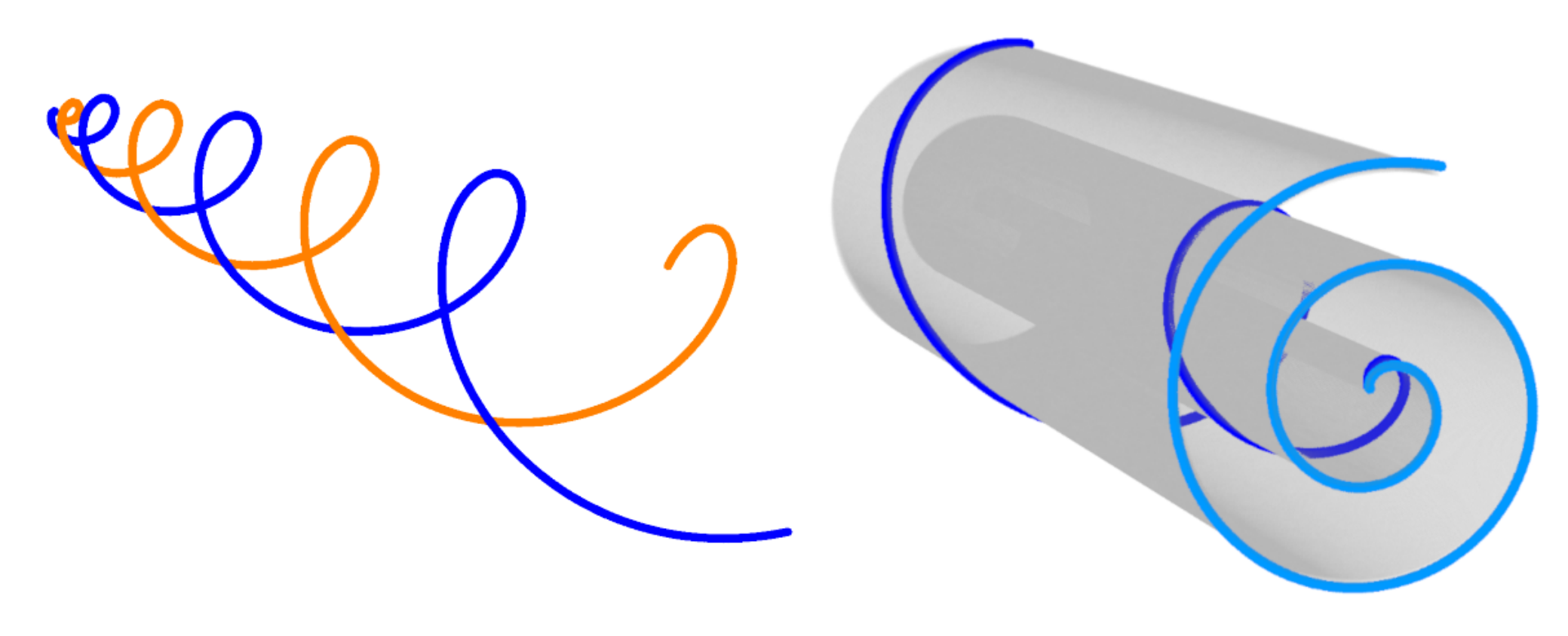}
\caption{A curve whose curvature and torsion are both equal to $c t^{-1/2}$ (blue) together with its congruent evolute (orange). The curves have constant slope and lie on generalized cylinders over the involutes of appropriate circles.}
\label{spirals}
\end{figure}

Figure \ref{spirals} shows one of the curves whose curvature and torsion satisfy $k = \tau = c t^{-1/2}$ for some $c > 0$, together with its congruent evolute. The curves have constant slope and are geodesics on generalized cylinders that have the involute of the circle with radius $\frac{1}{4\sqrt{2}c^2}$ as base curve (see \cite{St} for more details). Thus, if the base curve $\gamma$ is parametrized by
\begin{equation*}
\gamma(t) = (\gamma_1 (t), \gamma_2(t)) = \frac{1}{4 \sqrt{2} c^2} ( \cos z + z \sin z, \ \sin z - z \cos z ), \ \ \text{where} \ \ z= 2 c \sqrt{2t},
\end{equation*}
then the corresponding curve with congruent evolute is given by
\begin{equation*}
t \mapsto ( \gamma_1(t), \gamma_2(t), \frac{1}{\sqrt{2}}t).
\end{equation*}

It is tempting to apply Corollary \ref{cor:TotalCurvTors} twice: 
the  curvature density and the torsion density of the second evolute coincide with those of the initial curve (under the non-restrictive assumption that the torsion is positive).
This is equivalent to the fact that the Frenet frame of the second evolute is parallel to the Frenet frame of the initial curve (which also follows from the identity $(t_e, n_e, b_e) = (\pm b, \pm n, \pm t)$ that we mentioned earlier).
\smallskip

{\bf Remark.} Two simultaneously parametrized curves with parallel tangents at the corresponding points are called \emph{Combescure transformations} of each other; for a detailed study see \cite{S}.
\smallskip

A curve is congruent to its second evolute if and only if an arclength parameter for the curve is also an arclength parameter for the second evolute.
This condition is equivalent to a complicated differential equation
\begin{equation}
\label{eqn:2ndEvolute}
\frac{1}{r\tau} \left(\frac{r'}{\tau}\right)' + \left( \frac{1}{\sigma\tau}\left(\frac{\sigma}{\tau}\right)' \right)' = 0.
\end{equation}

However, if the curvature $k$ is constant, then $\sigma = \frac{\tau}{k}$, and equation \eqref{eqn:2ndEvolute} is satisfied. It follows that a space curve of constant curvature is congruent to its second  evolute. See \cite{FT} for curves that are homothetic to their second evolutes.

\paragraph{Involutes.}

By definition, the curve $\xi$ is an  involute of a curve $e$ if $e$ is the evolute of $\xi$. 
A generic curve has a two-parameter family of  involutes: they are the curves orthogonal to the family of osculating planes of the given curve. 

Let $S$ be the tangent developable of $e$, and let $H$ be the osculating plane of $e$ at some point. Let us roll the plane $H$ along the surface $S$ without slipping and twisting. The instantaneous motion of the plane $H$ is a rotation about a line in $H$. Therefore the velocity of each point of the plane $H$ is perpendicular to $H$, and the involutes of $e$ are the trajectories of the points of $H$ in the process of  rolling.

This description shows that the involutes of a curve are equidistant: the distances between the points of the rolling plane do not change. It also shows that each involute has a cusp every time it reaches $S$.

Assume that $e$ is a generic smooth closed curve with non-vanishing torsion. After rolling the osculating plane $H$ all the way around $e$, this plane returns to the original position, but its points do not necessarily return to their initial positions: a self-map of $H$ arises, which we call the \textit{monodromy}. 

If we orient the osculating planes of a curve by its binormals, then the monodromy is an isometry of $H$ that preserves this orientation.  A generic orientation-preserving isometry of the plane is a rotation about some point, and this point is its only fixed point. We conclude that a generic smooth closed curve $e$ has a unique closed involute. 

This is in contrast with the planar case: for the involute of a closed plane curve to be closed, the curve must have zero alternating perimeter (the sign changes after every cusp), and if this alternating perimeter vanishes, then all involutes are closed; see \cite{GTT} or \cite[Chapter~10]{FTb}.

Consider the trace $e_H$ of the curve $e$ in the plane $H$ as this plane rolls along $e$. This trace is a planar development of the curve $e$: its curvature, as a function of the arclength, is the same as that of $e$. (Imagine that a curve is made of wire which is hard to bend but easy to twist. Then one can flatten this curve without changing its curvature.)

While the curve $e$ is closed, $e_H$ in general is not. The end points of $e_H$ and the tangent directions therein
are two contact elements in the plane $H$. 

\begin{proposition}
The monodromy is the orientation preserving isometry of $H$ which sends the initial contact element of the curve $e_H$ to its terminal contact element. 
\end{proposition}

{\bf Proof.} When one surface is rolled, without slipping and twisting, along a curve on another surface, this curve and its trace on the rolling surface have the same curvature at the respective points. This describes the relation between the curves $e$ and $e_H$.

An instantaneous displacement of the plane $H$ when it is rolled along the curve $e$ is the parallel translation distance $dt$ in the tangent direction $e'(t)$, combined with the rotation through the angle $k(t) dt$. Integrated along $e$, this defines the motion described in the statement.
\proofend

\begin{corollary}
The rotation angle of the monodromy is equal to the total curvature of the evolute.
\end{corollary}

Assume that a curve has a closed involute. When are all involutes also closed, that is, when is the monodromy the identity map?

\begin{corollary}
If a curve has at least one closed involute, and the total torsion of this involute is an integer multiple of $2\pi$, then all involutes are closed.
\end{corollary}

This follows from Corollary \ref{cor:TotalCurvTors}: the total curvature of the evolute is equal to the total absolute torsion of the curve.
The torsion of the involute does not change its sign provided the curve has no cusps.

An example of a closed space curve whose planar development is also closed is a curve of constant curvature 1 and length $2k\pi$ (the development is a $k$ times traversed circle).

\paragraph{Osculating circles.} As we mentioned earlier, the osculating circles of a plane curve with monotonic curvature are nested. What about the osculating circles of a space curve? In particular, can the osculating circles at neighboring points be linked?

\begin{figure}[ht]
\centering
\includegraphics[width=.9\textwidth]{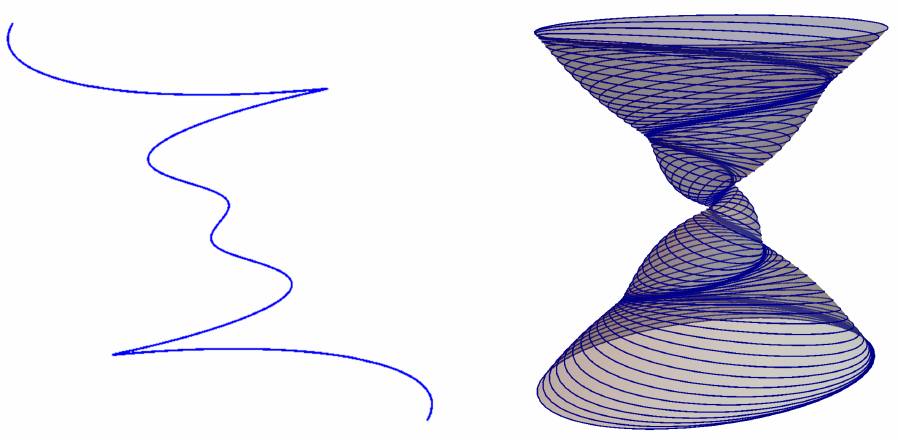}
\caption{Osculating circles of a space curve.}
\label{fig:OscCircles}
\end{figure}

The osculating circles of the curve in Figure \ref{fig:OscCircles} seem to be unlinked.
The next proposition confirms that, locally, this is always the case.

\begin{proposition} \label{link}
Let $\xi$ be a curve with non-vanishing torsion. Then its osculating circle at every point is disjoint from its osculating planes at sufficiently close points. In particular, the osculating circles at close points are not linked.
\end{proposition}

{\bf Proof.}  We view the curve as the envelope of the 1-parameter family of planes $z=a(t)x+b(t)y+c(t)$. We  assume that, for $t=0$, the plane is $z=0$, that $\xi(0)=(0,0,0)$, and the tangent line to $\xi$ at the origin is the $x$-axis. This means that
$$
a(0)=b(0)=c(0)=a'(0)=c'(0)=c''(0)=0.
$$
Thus
$$
a(t)=at^2+O(t^3),\ b(t)=bt+O(t^2),\ c(t)=ct^3+O(t^4).
$$
We will ignore the ``big O" terms in what follows: this will not affect the result, but will make the formulas less awkward.

We find the equation of the curve using the middle equation (2):
$$
\xi(t)=\left(-\frac{3ct}{a},\frac{3ct^2}{b},ct^3  \right).
$$
The curvature of this curve at the origin equals $\frac{2a^2}{3bc}$, and the radius of curvature is $\frac{3bc}{2a^2}$.

The intersection of the osculating planes of the curve $\xi(t)$, that is, the planes $z=a(t)x+b(t)y+c(t)$, with the plane $z=0$ are the lines $a(t)x+b(t)y+c(t)=0$ (when $t=0$, it is the $x$-axis), that is, $at^2x+bty+ct^3=0$.

 The envelope of this family of lines is the parabola
 $$
 x=\-\frac{2ct}{a},\ y=\frac{ct^2}{b}.
 $$
 The curvature of this curve at the origin equals $\frac{a^2}{2bc}$, and its radius of curvature  $\frac{2bc}{a^2} > \frac{3bc}{2a^2}$. 
 
 \begin{figure}[hbtp]
\centering
\includegraphics[height=3in]{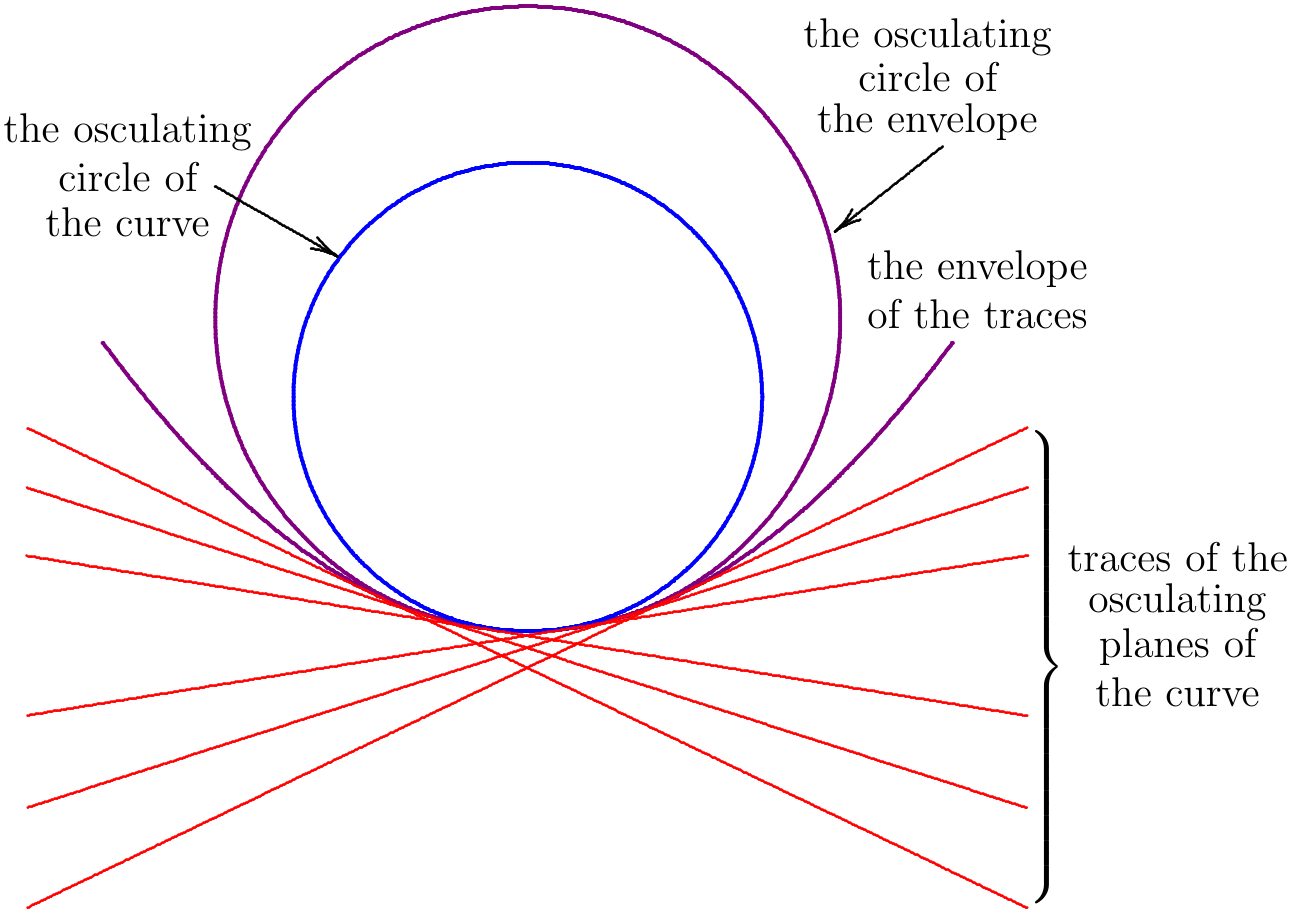}
\caption{To proof of Proposition \ref{link}.}
\label{para}
\end{figure}
 
 It follows that the osculating circle of the curve $\xi$ at the origin lies inside the parabola. Therefore the tangent lines to the parabola, that is, the intersections of the osculating planes of $\xi(t)$ with its osculating plane at $\xi(0)$, are disjoint from the osculating circle at $\xi(0)$.
This implies the result, see Figure \ref{para}. 
\proofend

\section{Pseudo-evolutes}\label{pseudo}

%\subsection{Definition.} 
\paragraph{Definition and first properties.}

A generic space curve $\xi$ determines a one-parameter family of rectifying planes. This  is the family of tangent planes to a certain developable surface $S$, called {\it rectifying developable}, and this surface has a regression edge. We call this regression edge the {\it pseudo-evolute} of the initial curve. In other words, the pseudo-evolute of a curve $\xi$ is a curve whose osculating planes are the rectifying planes of $\xi$.

\begin{proposition} \label{geodesic}The curve $\xi$ is a geodesic on $S$.\end{proposition}

{\bf Proof.} According to one of many description of geodesics on a surface in space, it is a curve whose principal normals are normal to the surface. For $\xi$ on $S$ this clearly holds: the principal normals of $\xi$ are perpendicular to the rectifying planes of $\xi$, that is, to the tangent planes of $S$.
\proofend

Informally, this means that a piece of paper with a straight line drawn on it may be, in a unique way, attached (without crumblings and foldings) to a given space curve in such a way that the line follows the curve (see \cite{FT1}).
This surface is unique because the normals to this surface are determined: they are the normals of $\xi$.
This gives another description of the pseudo-evolute of a (generic) curve $\xi$: it is the regression edge of the unique developable surface $S$, which contains $\xi$ as a geodesic.

\paragraph{Equation of the pseudo-evolute.}

The equation of the rectifying plane to the curve $\xi$ at the point $\xi(t)$ is
$$
(P-\xi(t))\cdot {\bf n}(t)=0.
$$
The first and second derivatives of this equations with respect to $t$ are
$$
\begin{array} {rll} &\hskip-6pt(P-\xi(t))\cdot(-k{\bf t}(t)+\tau{\bf b}(t))&\hskip-6pt=0\\ k\, +&\hskip-6pt(P-\xi(t))\cdot(-k'{\bf t}+\tau'{\bf b}-(k^2+\tau^2){\bf n})&\hskip-6pt=0.\end{array}
$$
These three equations form the system that has the solution
$$
(P-\xi(t))\cdot {\bf t}(t)=\frac{k\tau}{k'\tau-k\tau'},\ (P-\xi(t))\cdot {\bf b}(t)=\frac{k^2}{k'\tau-k\tau'}.
$$
We obtain the equation of the pseudo-evolute $\varepsilon(t)$ of the curve $\xi(t)$: 

\begin{proposition} \label{equation-pseudo} 
The pseudo-evolute $\varepsilon$ of the curve $\xi$ has the equation
$$
\varepsilon=\xi+\frac k{k'\tau-k\tau'}(\tau{\bf t}+k{\bf b}).
$$
\end{proposition}

\paragraph{Escapes to infinity and cusps.}

\begin{proposition} \label{infcusp}
Let $\xi$ be a generic curve with curvature $k$ and torsion $\tau$.

$(1)$ The pseudo-evolute of $\xi$ escapes to infinity when the first derivative of the ratio $\displaystyle\frac{\tau}k$ vanishes.

$(2)$ The pseudo-evolute of $\xi$ has cusps when the second derivative of $\displaystyle\frac{\tau}k$ vanishes.
\end{proposition}

\begin{figure}[ht]
\centering
\includegraphics[width=.52\textwidth]{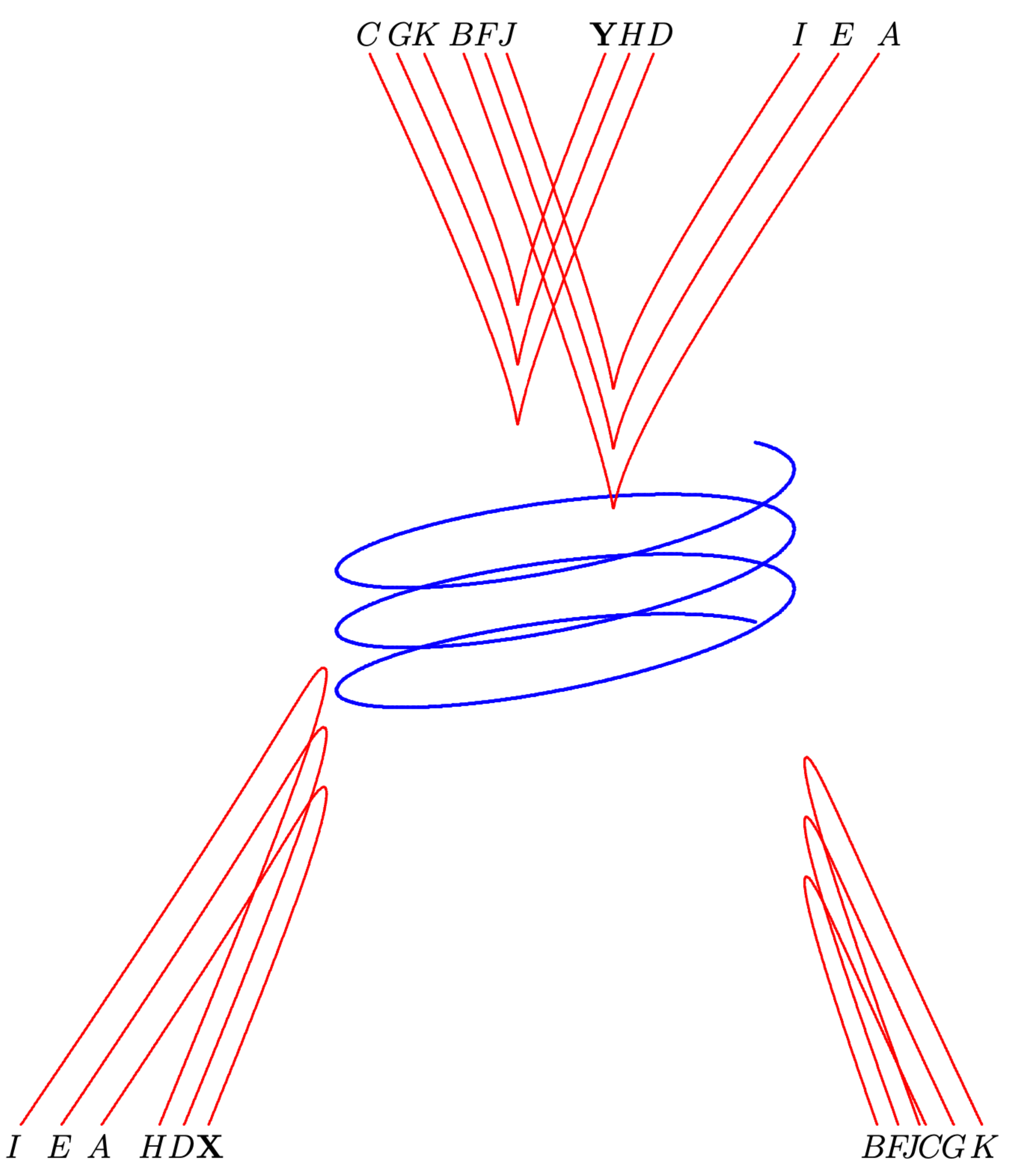}\ 
\includegraphics[width=.45\textwidth]{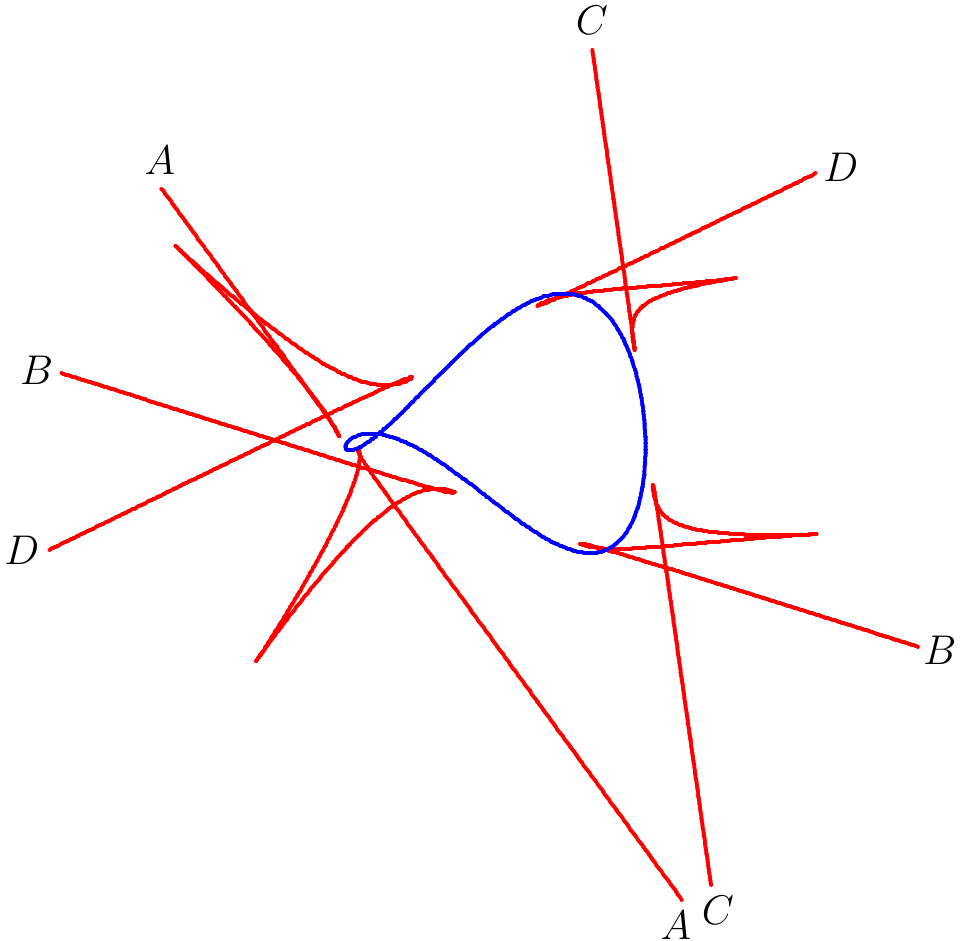}
\caption{Left: the pseudo-evolute of an elliptical helix. Right: the pseudo-evolute of a closed smooth curve. It has 12 cusps and escapes 4 times to infinity.}
\label{pevol}
\end{figure}

{\bf Proof.} Part (1) follows from Proposition \ref{equation-pseudo}: $\varepsilon$ escapes to infinity, when $k'\tau-k\tau'=0$.

Part (2) requires some computation also based on the formula from Proposition \ref{equation-pseudo}. From this formula: 
\begin{equation} \label{eq:old4}
\varepsilon'={\bf t}+\left(\frac k{ k'\tau-k\tau'}\right)'(\tau{\bf t}+k{\bf b})+\frac k{ k'\tau-k\tau'}(\tau{\bf t}+k{\bf b})'.\
\end{equation}
By Frenet's formulas, $(\tau{\bf t}+k{\bf b})'=(\tau'{\bf t}+ k'{\bf b})$, so the right hand side of (\ref{eq:old4}) is a linear combination of $\bf t$ and $\bf b$. The coefficient of $\bf t$ is
$$
\begin{array} {c} \displaystyle{1+\frac{( k'( k'\tau-k\tau')-k( k'\tau-k\tau')'\tau}{( k'\tau-k\tau')^2}+\frac{k\tau'}{ k'\tau-k\tau'}}\\ = \displaystyle{\frac{( k'\tau-k\tau')+ k'\tau+k\tau'}{ k'\tau-k\tau'}-\frac{k\tau( k'\tau-k\tau')'}{( k'\tau-k\tau')^2}=\left(\frac{2 k'}{ k'\tau-k\tau'}-\frac{k( k'\tau-k\tau')'}{( k'\tau-k\tau')^2}\right)\tau},\end{array}
$$
while the coefficient of $\bf b$ is 
$$
\frac{( k'( k'\tau-k\tau')-k( k'\tau-k\tau')'k}{( k'\tau-k\tau')^2}+\frac{k k'}{ k'\tau-k\tau'}=k\left(\frac{2 k'}{ k'\tau-k\tau'}-\frac{k( k'\tau-k\tau')'}{( k'\tau-k\tau')^2}\right).
$$
Hence, $\varepsilon'=0$ if $\displaystyle{\frac{2 k'}{ k'\tau-k\tau'}-\frac{k( k'\tau-k\tau')'}{( k'\tau-k\tau')^2}}=0$, and the last expression is $\displaystyle{\frac{k^3}{( k'\tau-k\tau')^2}\left(\frac\tau k\right)''}$.
\proofend

Thus the derivative $\displaystyle{\left(\frac\tau k\right)'}$ plays, for pseudo-evolutes, the role similar to that of the curvature for evolutes of planar curves. 

The function $\displaystyle{\frac\tau k}$ has the following geometrical meaning. The {\it tangent indicatrix} of a space curve $\xi(t)$ is the curve $\xi'(t)$ on the unit sphere (recall that $\xi$ is arclength-parametrized). The geodesic curvature of the tangent indicatrix equals $\displaystyle{\frac\tau k}$.

Proposition \ref{infcusp} has the following consequence (\cite[Section 2-4, Problem 5]{St}):

\begin{corollary}
The rectifying developable is a cone if and only if $\displaystyle\frac{\tau}{k}$ is a linear function, and it is a cylinder if and only if $\displaystyle\frac{\tau}{k}$ is a constant.
\end{corollary}

\paragraph{Examples.}

Figure \ref{pevol} shows two examples: the pseudo-evolutes of the same elliptical helix as in Figure \ref{ex23}, and of the curve $x=\cos t, y=\sin t, z=\frac{1}{2}\sin 2t$. 

\begin{figure}[ht]
\centering
\includegraphics[width=.65\textwidth]{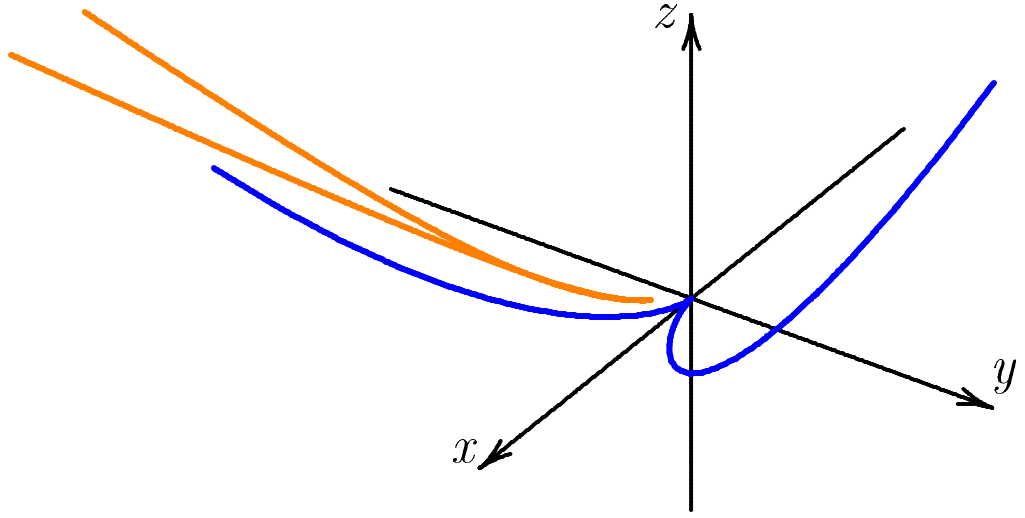}
\caption{The pseudo-evolute of a curve with a cusp.}
\label{pevolcusp}
\end{figure}

\paragraph{Pseudo-evolute of a curve with a cusp.}
Our definition works in this case, but we are  puzzled by the geometry of the result, see Figure \ref{pevolcusp}.

We consider the standard example of the curve with a cusp: $x=t^2,y=t^3,z=t^4$. It is true that the principal normal of a curve (at a regular point) has the directiion of $\bf t\times(t\times\dot t)$ which, for our curve, is
$$
(18t^4+64t^6,48t^7-12t^3,-32t^4-36t^6)=2t^3(9t+32t^3,24t^4-6,-16t-18t^3).
$$
Thus, the rectifying plane of our curve at a point with the parameter value $t\neq0$ has the equatiion
$$
\begin{array} {rl} F(x,y,z,t)=(9t+32t^3)(x-t^2)+(24t^4-6)(y-t^3)-(16t+18t^3)(z-t^2)&\\ =(9t+32t^3)x+(24t^4-6)y-(16t+18t^3)z-(6t^7+16t^5+3t^3)=0.&\end{array}
$$
This equation makes sense for $t=0$ (it becomes $y=0$), so we have a family of rectifying planes parametrized by $t\in\mathbb R$. 

To find the envelope of this family, we consider the system $F=F_t=F_{tt}=0$, and we %have a temptation to 
cancel the equation $F_{tt}=0$ by $t$. After doing that, we solve the system using the Cramer rule and obtain:
\begin{equation} \label{monster}
\begin{array} {rl} x=&\displaystyle\frac{4(324t^{10}+1440t^8+1335t^6+720t^4+320t^2+18)}{525(1-4t^4)},\\ y=&-\displaystyle\frac{8t^5(9t^2+8)}{9(1-4t^4)},\\z=&\displaystyle\frac{460t^{10}+7880t^8+11520t^6+8490t^4+1440t^2+81}{1050(1-4t^4)}.\end{array}
\end{equation}
These are parametric equations of a curve that can be considered as the pseudo-evolute of the curve $x=t^2,y=t^3,z=t^4$. 

The curve (\ref{monster}) has a singularity at the point $\displaystyle\frac3{350}(16,9)$, but this singularity is a ``degenerate cusp":  it is easy to find that $$16z-9x=\displaystyle\frac{4(36t^{10}+16t^8+63t^6+60t^4)}{3(1-4t^4)},$$which shows that the ``type of the cusp" is $(t^2,t^4,t^5)$. How to relate this to our geometric interpretation of pseudo-evolutes remains unclear to us.

%Notice also that if we do not cancel the equation $F_{tt}=0$ by $t$, then our system gains an additional solution $16z-9x=y=t=90$, so maybe we need to count the line $16z-0x=y=0$ as a part of the pseudo-evolute, which all corresponds to the parameter value $t=0$. This, however, does not have much to do with the understanding of the geometry of the pseudo-evolute of the curve $x=t^2,y=t^3,z=t^4$.

\paragraph{Pseudo-involutes.}
Every curve has a two-parameter family of pseudo-involutes, namely the curves whose pseudo-evolute is the initial curve; indeed, a curve $\eta$ is the regression edge of a developable surface, the union of tangent lines of $\eta$, and the geodesics of this surface are pseudo-involutes of $\eta$. 

Pseudo-involutes may have cusps: it happens every time when the involute reaches $\eta$. (This makes  pseudo-involutes similar to involutes of planar curves, but there is a big difference: pseudo-involutes of $\eta$ do not need to be {\it perpendicular} to $\eta$ at the cuspidal points.)

A pseudo-involute may be a smooth closed curve -- see Figure \ref{pevol}, right. Then it is a closed geodesic on a developable surface, and its equidistant curves on this surface are also closed geodesics. Hence such a curve is included in a 1-parameter family of smooth closed pseudo-involutes.

\paragraph{An annoying question.}
We see that pseudo-evolutes have many properties similar to those of evolutes. Still their geometric meaning remains enigmatic. Evolutes are the loci of centers of osculating spheres. And pseudo-evolutes are the loci -- of what?\footnote{Spoiler alert: the answer will be given at the very end of the article.}

\section{Monge evolutes} \label{sect:Monge}

\paragraph{Monge involutes.} Attach one end of a non-stretchable string to a space curve $\eta$, pull the string tight, and wrap it around $\eta$. The velocity of the free end of the string is always orthogonal to the string (it is non-stretchable). The trajectory of this free end is the curve $\xi$, a Monge involute of $\eta$. Changing the length of the string yields a 1-parameter family of these involutes.

Differentiating the formula
\begin{equation}
\label{eqn:StringInvolute}
\xi(t) = \eta(t) + (\ell-t) {\bf t}_\eta,
\end{equation}
where $\ell$ is the length of the string attached at point $\eta(0)$,
yields
\begin{equation}
\label{eqn:XiPrime}
\xi' = (\ell-t) k_\eta {\bf n}_\eta.
\end{equation}
It follows that if the curvature of $\eta$ never vanishes, then the cusps of $\xi$ lie on $\eta$, which happens when the free end of the string lands on the curve.

\paragraph{Monge evolutes.} If $\xi$ is a Monge involute of $\eta$ then, by definition, $\eta$ is a Monge evolute of $\xi$. 
Does every curve has such an evolute, and if so, how many?

If $\eta$ is an evolute of $\xi$, then one has
\[
\eta = \xi + y {\bf n} + z {\bf b}.
\]
Differentiate with respect to the arclength parameter of $\xi$ and apply the Frenet formulas:
\[
\eta' = (1 - k y) {\bf t} + (y' - \tau z){\bf n} + (z' + \tau y){\bf b}.
\]
By assumption, $\eta'$ must be parallel to $\eta - \xi = y {\bf n} + z {\bf b}$, which leads to a system of equations
\[
1 - ky = 0, \quad
\frac{y' - \tau z}{y} = \frac{z' + \tau y}{z}.
\]
This system has a one-parameter family of solutions
\[
y = r, \quad z = -r \, \tan \int \tau\, dt.
\]

This leads to
\begin{proposition}
\label{prp:StringEvolute}
The Monge evolutes of an arclength-parametrized curve $\xi$ are given by
\[
\eta = \xi + r{\bf n} - r \, \tan\alpha\, {\bf b}, \quad \alpha' = \tau.
\]
\end{proposition}

In particular, spatial Monge evolutes of a plane curve are geodesics on the cylinder over the plane evolute of this curve.

\begin{corollary}
The corresponding points of two different evolutes are seen from the corresponding point of the curve under a constant angle.
\end{corollary}

\paragraph{Singularities of Monge evolutes.}
Proposition \ref{prp:StringEvolute} implies that
\[
\|\eta - \xi\| = \frac{1}{k|\cos\alpha|}.
\]
On the other hand, due to \eqref{eqn:StringInvolute}, one has
\[
\|\eta - \xi\| = \int_{t_0}^t \|\eta'\|\, dt.
\]
Therefore
\begin{equation*}
\label{eqn:NormEtaPrime}
\|\eta'\| = \left( \frac{1}{k|\cos\alpha|} \right)'.
\end{equation*}

It follows  that the cusps of a Monge evolute correspond to the critical points of the function $k \cos\alpha$.
This generalizes a property of plane evolutes: their cusps correspond to the vertices of the curve.

Proposition \ref{prp:StringEvolute} also allows to determine when an evolute escapes to infinity: this happens if either $k = 0$ or $\alpha = \frac{\pi}2 + k\pi$.

\paragraph{Closed Monge evolutes and Monge involutes.}
Proposition \ref{prp:StringEvolute} has the following consequence.

\begin{proposition}
Monge evolutes of a closed space curve are closed if and only if the total torsion of the curve is an integer multiple of $\pi$.
In particular, Monge evolutes of a centrally symmetric curve are closed.
\end{proposition}

The last statement is due to the obvious fact that a central symmetry reverses the sign of the torsion.

\begin{figure}[ht]
\centering
\includegraphics[width=.9\textwidth]{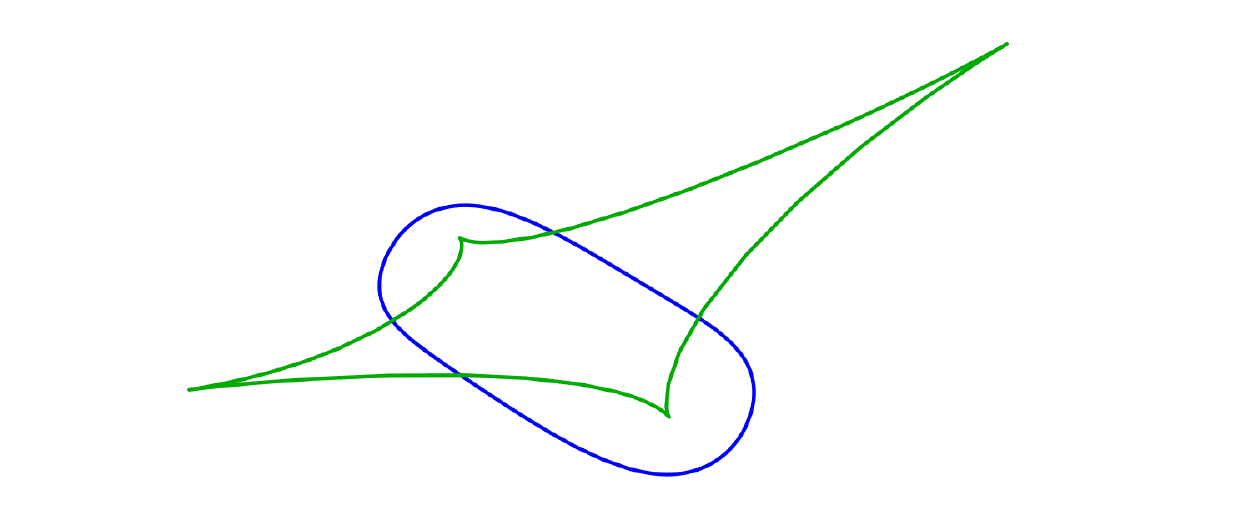}
\caption{A Monge evolute (green) of a centrally symmetric curve.}
\end{figure}

As to the closed Monge involutes, the situation is exactly the same as in the plane:

\begin{proposition}
Monge involutes of a closed space curve are closed if and only if the curve has zero length (the length element changes sign after each cusp).
\end{proposition}

\paragraph{Interrelations between the evolutes.} As promised, we describe how the three kinds of evolutes interact with each other.

The line $\{\xi + r{\bf n} + \lambda{\bf b} \mid \lambda \in {\mathbb R}\}$ is called the \emph{polar line} of the curve $\xi$ at the corresponding point.
This is the line that goes through the center $\xi + r{\bf n}$ of the osculating circle and is orthogonal to the osculating plane of $\xi$.
The center of the osculating sphere lies on the polar line.

Recall the formula $e' = \sigma{\bf b}$. It implies that the polar line is tangent to the evolute of $\xi$.
It follows that the polar lines are the rulings of the normal developable  of the curve $\xi$.

Proposition \ref{prp:StringEvolute} implies that each point of a Monge evolute lies on some polar line,  and every point of a polar line lies on a unique Monge evolute.
It follows that the  normal developable surface  of the curve $\xi$ is foliated by its Monge evolutes.

A Monge evolute of $\xi$ meets the evolute of $\xi$ at the points where $r' = r\tau \tan\alpha$.
The tangent to the Monge evolute at those points coincides with the binormal of $\xi$, that is, the Monge evolute is tangent to the evolute. 

As a result we have

\begin{proposition}
The evolute is the envelope of Monge evolutes.
\end{proposition}

\begin{figure}[ht]
\centering
\hspace*{-3.1cm}\includegraphics[scale=0.8]{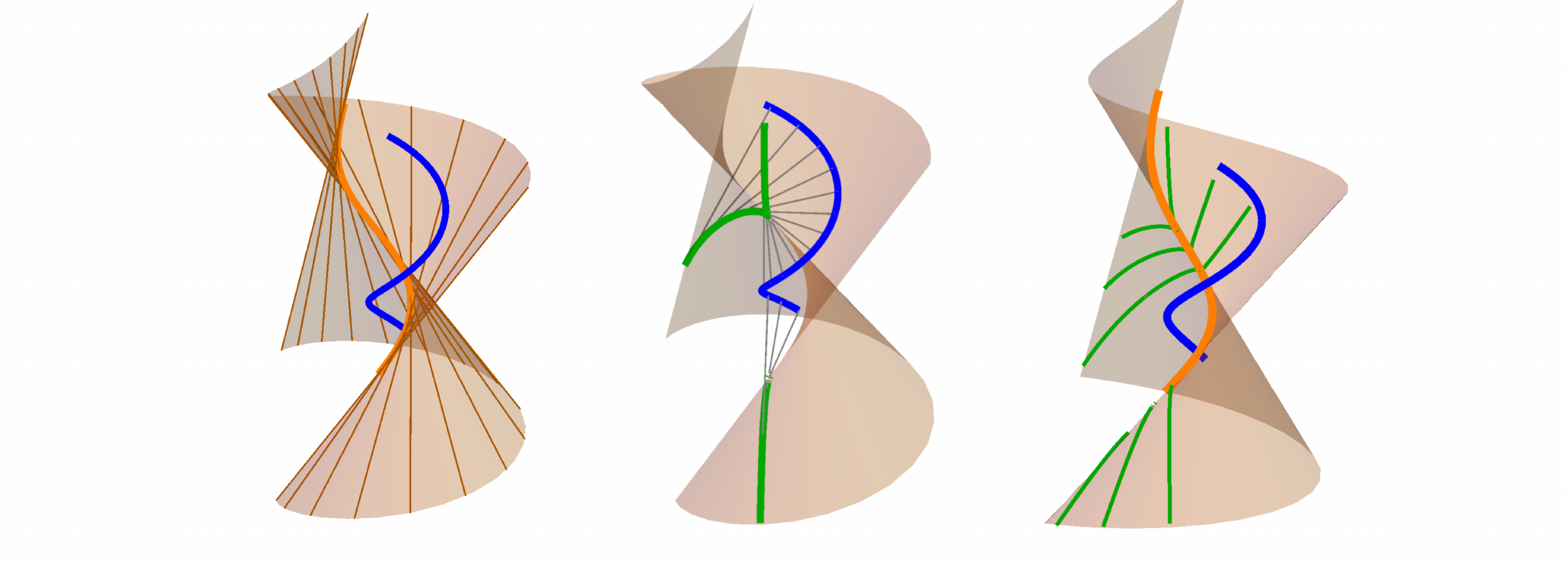}
\caption{The evolute (orange) and some Monge evolutes (green) of a  curve~$\xi$ (blue) together with its normal developable (orange surface). \textit{Left:} The polar lines are tangent to the evolute and lie on the normal developable of $\xi$. \textit{Middle:} The tangents of the Monge evolute intersect $\xi$ orthogonally. \textit{Right:} The Monge evolutes also foliate the normal developable of $\xi$.}
\end{figure}

We conclude with yet another relation between the different types of evolutes (\cite[Section 2-4, Exercise 3]{St}).

\begin{proposition}
The pseudo-evolute of any space curve is the  evolute of any of its Monge involutes.
\end{proposition}

{\bf Proof.}
Let $\eta$ be a space curve, and let $\xi$ be a Monge evolute of $\eta$.
One has to show that the normal planes of $\eta$ are the rectifying planes of $\xi$.
These planes are parallel because their normals are parallel; see Equation \eqref{eqn:XiPrime}.
And they coincide because $\xi-\eta$ is tangent to $\eta$ and normal to $\xi$, and thus is contained in both planes.
\proofend

In particular, this provides an answer to the above ``annoying question": the pseudo-evolute of a curve  is the locus of the osculating spheres of its Monge evolute.

\paragraph{Acknowledgements.} I. Izmestiev and M. Raffaelli were supported by the Austrian Science Fund (FWF), project F 77 (SFB "Advanced Computational Design"); S. Tabachnikov was supported by NSF grant  DMS-2005444.


\begin{thebibliography}{99}

\bibitem{B} W. Blaschke. {\it Vorlesungen \"uber Differentialgeometrie}. Springer-Verlag, Berlin, 1930.

%\bibitem{BL} W. Blaschke, K. Leichweiss, {\it Elementarte Differentialgeometrie}, Springer-Verlag, Berlin, 1973.

\bibitem{BJT} G. Bor, C. Jackman, S. Tabachnikov.
{\it Variations on the Tait-Kneser theorem.} 
Math. Intelligencer {\bf 43} (2021), no. 3, 8--14. 

\bibitem{CSW} G. Cairns, R. Sharpe, L. Webb,.
{\it Conformal invariants for curves and surfaces in three-dimensional space forms.} 
Rocky Mountain J. Math. {\bf 24} (1994), 933--959. 

\bibitem{DO} E. Demaine, J. O'Rourke. {\it Geometric folding algorithms. Linkages, origami, polyhedra.} Cambridge Univ. Press, Cambridge, 2007. 

\bibitem{E} L.  Eisenhart. {\it A treatise on the differential geometry of curves and surfaces.} Dover Publications, Inc., New York 1960.

\bibitem{F} D. Fuchs. {\it Evolutes and involutes of spatial curves}. Amer. Math. Monthly. {\bf120} (2013), 217--231.

\bibitem{FT1} D. Fuchs, S. Tabachnikov. {\it More on paperfolding.} Amer. Math. Monthly {\bf 106} (1999), 27--35. 

\bibitem{FT} D. Fuchs, S. Tabachnikov. {\it Iterating evolutes of spacial polygons and spacial curves}. Moscow Math. J. {\bf120}  (2017) , 667--689.

\bibitem{FTb} D. Fuchs, S. Tabachnikov. {\it  Mathematical omnibus. Thirty lectures on classic mathematics.} Amer. Math. Society, Providence, RI, 2007. 

\bibitem{GTT} E. Ghys, S. Tabachnikov, V. Timorin.
{\it Osculating curves: around the Tait-Kneser theorem.} 
Math. Intelligencer {\bf 35} (2013), no. 1, 61--66. 

\bibitem{La} J. Langer. {\it Recursion in curve geometry.} New York J. Math. {\bf 5} (1999), 25--51.

\bibitem{S} E. Salkowski {\it Zur Transformation von Raumkurven}. Math. Annalen {\bf 66}, (1909), 517--557.

\bibitem{St} D. Struik. {\it Lectures on classical differential geometry.} Reprint of the second edition. Dover Publications, Inc., New York, 1988.

\bibitem{UV} R. Uribe-Vargas. {\it On vertices, focal curvatures, and differential geometry of space curves}. Bull. Braz. Math. Soc. {\bf36}  (2005), 285--307.

\end{thebibliography}
\end{document}